\documentclass[11pt,reqno]{amsart}
\usepackage{amsmath}
\usepackage{amsthm}
\usepackage{graphicx}
\usepackage{amsfonts}
\usepackage{amssymb}
\usepackage{enumerate}
\usepackage[margin=1.1in]{geometry}
\numberwithin{equation}{section}

\newtheorem{theorem}{Theorem}[section]

\theoremstyle{definition}

\newtheorem{definition}[theorem]{Definition}
\newtheorem{remark}[theorem]{Remark}
\newtheorem*{assumption*}{Assumption}

\def\E{{\mathbb E}}
\def\R{{\mathbb R}}
\def\N{{\mathbb N}}
\def\FF{{\mathbb F}}
\def\PP{{\mathbb P}}
\def\EE{{\mathcal E}}

\def\P{{\mathcal P}}

\def\H{{H}}

\def\W{{\mathcal W}}

\def\F{{\mathcal F}}

\def\C{{\mathcal C}}

\title{On a strong form of propagation of chaos for McKean-Vlasov equations}
\author{Daniel Lacker}

\begin{document}

\begin{abstract}
This note shows how to considerably strengthen the usual mode of convergence of an $n$-particle system to its McKean-Vlasov limit, often known as propagation of chaos, when the volatility coefficient is nondegenerate and involves no interaction term. Notably, the empirical measure converges in a much stronger topology than weak convergence, and any fixed $k$ particles converge in total variation to their limit law as $n\rightarrow\infty$. This requires minimal continuity for the drift in both the space and measure variables. The proofs are purely probabilistic and rather short, relying on Girsanov's and Sanov's theorems. Along the way, some modest new existence and uniqueness results for McKean-Vlasov equations are derived.
\end{abstract}

\maketitle

\section{Introduction}

This note develops a simple but apparently new approach to analyzing McKean-Vlasov stochastic differential equations, of the form
\[
dX_t = b(t,X_t,\mu_t)dt + \sigma(t,X_t)dW_t, \quad \mu_t = \mathrm{Law}(X_t), \quad \forall t \ge 0,
\]
in which the drift is merely bounded and measurable, with fairly weak continuity requirements in the measure variable. The volatility $\sigma$ is nondegenerate and independent of the measure, and this enables a line of argument based on Girsanov's theorem which leads to a much stronger propagation of chaos result than usual, along with some new results on existence and uniqueness.

Propagation of chaos here refers to the convergence of the $n$-particle system, defined by the SDE
\begin{align*}
dX^{n,i}_t &= b\Big(t,X^{n,i}_t,\mu^n_t\Big)dt + \sigma(t,X^{n,i}_t)dW^i_t, \quad\quad\quad \mu^n_t = \frac{1}{n}\sum_{j=1}^n\delta_{X^{n,j}_t},
\end{align*}
to the solution law $\mu$ of the McKean-Vlasov equation. Precisely, propagation of chaos typically means that the empirical measures $\mu^n$ (say, on the path space) converge weakly in probability to the deterministic measure $\mu$, or equivalently that the law of $(X^{n,1},\ldots,X^{n,k})$ converges weakly to the product measure $\mu^{\otimes k}$ for any fixed $k$. In our context, we show that in fact the law of $(X^{n,1},\ldots,X^{n,k})$ converges \emph{in total variation} to $\mu^{\otimes k}$. Moreover, the sense in which $\mu^n$ converges in probability to $\mu$ can be strengthened; rather than working with the usual weak topology induced by duality with  bounded continuous test functions, we work with the stronger topology induced by duality with bounded measurable test functions.
In particular, our results will assume that $b(t,x,\mu)$ is continuous in $\mu$ in this stronger topology (or, for some results, in total variation) and merely measurable in $(t,x)$, and our coefficients may be path-dependent as well. 

McKean-Vlasov equations have been studied in a variety of contexts since the seminal work of McKean \cite{McKean1907}. Sznitman's monograph \cite{sznitman1991topics} is a classic introduction, and G\"artner's results \cite{gartner1988mckean} remain among the most general on existence, uniqueness, and propagation of chaos results for models with (weakly) continuous coefficients.

More recently, interacting diffusion models of this form have enjoyed something of a renaissance, due in part (but certainly not entirely) to new applications in \emph{mean field game theory} \cite{lasrylions}, and this is one impetus for revisiting these classical questions here. The McKean-Vlasov equations arising in mean field game theory can involve feedback controls obtained via Nash equilibrium problems. Regularity for these controls can be hard to come by, and this motivates a better understanding of somewhat more pathological dynamics. For instance, the recent work of \cite{campi-fischer} on mean field games with absorbing states naturally gives rise to McKean-Vlasov systems with discontinuous and path-dependent coefficients.

Several authors have studied McKean-Vlasov systems with various kinds of discontinuities arising in a variety of concrete applications. Noteworthy classes of examples include interactions based on ranks \cite{shkolnikov2012large,jourdain2013propagation} and quantiles \cite{crisan2014conditional,kolokoltsov2013nonlinear},  to which our results apply in certain cases. One such example given in Section \ref{se:example}, where we show that the particle approximation of Burgers' equation given in \cite{bossy-talay2,jourdain1997diffusions} holds in a stronger sense.

While several  papers have studied McKean-Vlasov equations with discontinuities, the coefficients are often \emph{continuous enough}, in the sense that the set of discontinuities has measure zero with respect to any candidate solution (see, e.g., \cite{chiang1994mckean}). In such a situation one can still apply the usual weak convergence arguments, which are not available for the general discontinuities in $x$ we allow, more in the spirit of \cite{jourdain1997diffusions}. We lastly mention the interesting recent works \cite{mishura2016existence,raynal2015strong} that deal  with similarly irregular coefficients but less general interaction terms, with no results on propagation of chaos. While our existence and uniqueness results differ from those mentioned above, the main novelty of this work is the strong propagation of chaos result, Theorem \ref{theorem-propagation}.

Section \ref{se:mainresults} below states the main results, and proofs are given in Sections \ref{se:proofs-existenceuniqueness} and \ref{se:proof-propagation}.
It is worth stressing that all of the proofs are purely probabilistic.

\section{Main results} \label{se:mainresults}

\subsection{Notation and topologies}
Let $E$ be a Polish space. For a signed Borel measure $\gamma$ on $E$, define the total variation norm
\[
\|\gamma\|_{\mathrm{TV}} := \sup\left\{\int_E f\,d\gamma : f : E \rightarrow \R \text{ measurable}, \ |f(x)| \le 1 \ \forall x \in E\right\}.
\]
Let $\P(E)$ denote the set of Borel probability measures on $E$. 
For $\mu,\nu \in \P(E)$, define the relative entropy
\begin{align}
\H(\mu | \nu) = \int_E \frac{d\mu}{d\nu}\log \frac{d\mu}{d\nu}\,d\nu, \ \  \text{ if } \mu \ll \nu, \quad\quad \H(\mu | \nu) =\infty \ \  \text{ otherwise.} \label{entropydef}
\end{align}
Let $B(E)$ denote the set of bounded measurable real-valued functions on $E$. Define $\tau(E)$ to be the coarsest topology on $\P(E)$ such that the map $\mu \mapsto \int_E\phi \, d\mu$ is continuous for each $\phi \in B(E)$. This topology is somewhat well known in large deviations literature as the \emph{$\tau$-topology}. Notably, $(\P(E),\tau(E))$ is not separable or metrizable.

The map $E^n \ni (x_1,\ldots,x_n) \mapsto \frac{1}{n}\sum_{j=1}^n\delta_{x_j} \in \P(E)$ need not be measurable with respect to the Borel $\sigma$-field of $(\P(E),\tau(E))$, and we will need to work with a smaller $\sigma$-field for which we recover this measurability. Define $\EE(\P(E))$ to be the smallest $\sigma$-field on $\P(E)$ such that the map $\mu \mapsto \int_E\phi \, d\mu$ is measurable for each $\phi \in B(E)$. It is well known that $\EE(\P(E))$ coincides with the Borel $\sigma$-field on $\P(E)$ generated by the topology of weak convergence \cite[Corollary 7.29.1]{bertsekasshreve}.

\subsection{The McKean-Vlasov equation} \label{se:existenceuniqueness}
Fix a time horizon $T > 0$ and a dimension $d \in \N$. Let $\C = C([0,T];\R^d)$ denote the path space, endowed with the supremum norm. We will be interested in McKean-Vlasov equations of the form
\begin{align}
dX_t = b(t,X,\mu)dt + \sigma(t,X)dW_t, \quad X_0 \sim \lambda_0, \quad \mu = \mathrm{Law}(X), \label{SDE}
\end{align}
stated more precisely in Definition \ref{def:weaksolution} below.
The data of the problem are coefficients
\begin{align*}
b &: [0,T] \times \C \times \P(\C) \rightarrow \R^d, \\
\sigma &: [0,T] \times \C \rightarrow \R^{d \times d},
\end{align*}
and an initial law $\lambda_0 \in \P(\R^d)$.

For $\mu \in \P(\C)$ and $t \in [0,T]$, let $\mu^t \in \P(\C)$ denote the law of the process stopped at time $t$, defined as the image of $\mu$ through the map $\C \ni x \mapsto x_{\cdot \wedge t} \in \C$.
At various points in the sequel, we will refer to the following assumptions:
\begin{itemize}
\item[($\mathcal{E}$)] $b$ is jointly measurable with respect to $\mathrm{Borel}([0,T]) \otimes \mathrm{Borel}(\C) \otimes \mathcal{E}(\P(\C))$, and $\sigma$ is jointly Borel-measurable. In addition, the coefficients are progressive in the sense that $\sigma(t,x) = \sigma(t,x_{\cdot \wedge t})$ and $b(t,x,\mu) = b(t,x_{\cdot \wedge t}, \mu^t)$ for every $(t,x,\mu)$.
\item[(A)] For each $(t,x)$ the matrix $\sigma(t,x)$ is invertible, and there exists $c > 0$ such that $|\sigma^{-1}b| \le c$. Moreover, there exists  a unique strong solution to the driftless SDE,
\[
dX_t = \sigma(t,X)dW_t, \quad X_0 \sim \lambda_0.
\]
\item[(B1)] $b(t,x,\cdot)$ is $\|\cdot\|_{\mathrm{TV}}$-Lipschitz, uniformly in $(t,x)$. More precisely, there exists $\kappa > 0$ such  that, for all $(t,x,\mu)$,
\[
|b(t,x,\mu) - b(t,x,\nu)| \le \kappa\|\mu^t - \nu^t\|_{\mathrm{TV}}.
\]
\item[(B2)] For each $\mu \in \P(\C)$, the following function is sequentially $\tau(\C)$-continuous at $\mu$:
\[
\P(\C) \ni \nu \mapsto \int_\C \int_0^T\left|\sigma^{-1}b(t,x,\nu) - \sigma^{-1}b(t,x,\mu)\right|^2\,dt\,\mu(dx).
\]
\end{itemize}

\begin{remark}
If one is careful about integrability, the assumptions can undoubtedly be relaxed to cover unbounded coefficients and stronger topologies for the continuity of $b(t,x,\mu)$ in $\mu$. We prefer to avoid obscuring the main line of argument with such generalities.
\end{remark}

\begin{definition} \label{def:weaksolution}
We say $\mu \in \P(\C)$ is a weak solution of \eqref{SDE} if there exists a filtered probability space $(\Omega,\F,\mathbb{F},\PP)$ supporting a progressively measurable $d$-dimensional process $X$, a $d$-dimensional  $\mathbb{F}$-Wiener process $W$, and an $\F_0$-measurable random vector $\xi$ with law $\lambda_0$, such that $\PP \circ X^{-1} = \mu$ and
\[
X_t = \xi + \int_0^tb(s,X,\mu)ds + \int_0^t\sigma(s,X)dW_s, \ \ \ t \in [0,T].
\]
\end{definition}

\begin{theorem} \label{existence-uniqueness}
Assume ($\mathcal{E}$), (A), and (B1) hold. Then \eqref{SDE} has a unique weak solution.
\end{theorem}

\begin{theorem} \label{existence}
Assume ($\mathcal{E}$), (A), and (B2) hold. Then \eqref{SDE} has a weak solution.
\end{theorem}

The closest result to Theorem \ref{existence} that we know of seems to come from the paper \cite{carmona-lacker}, from which we borrow the proof idea. A nearly identical form of Theorem \ref{existence-uniqueness} was given in \cite[Theorem 2.2]{jourdain1997diffusions} and \cite[Theorem C.1]{campi-fischer}, though our proof seems to be much simpler.

\subsection{Propagation of chaos} \label{se:propagation}

For $n \in \N$, let $(X^{n,1},\ldots,X^{n,n})$ denote a weak solution on some filtered probability space $(\Omega,\F,\FF,\PP)$ of the SDE system
\begin{align*}
dX^{n,i}_t &= b(t,X^{n,i},\widehat\mu^n)dt + \sigma(t,X^{n,i} )dW^i_t, \quad X^{n,i}_0 = \xi_i, \\
\widehat\mu^n &= \frac{1}{n}\sum_{k=1}^n\delta_{X^{n,k}},
\end{align*}
where $W^1,\ldots,W^n$ are independent $d$-dimensional $\FF$-Wiener processes, and $\xi_1,\ldots,\xi_n$ are i.i.d.\ and $\F_0$-measurable with law $\lambda_0$. Under assumptions ($\mathcal{E}$) and (A), a standard argument by Girsanov's theorem guarantees the existence and uniqueness in law for this SDE system.

\begin{theorem} \label{theorem-propagation}
Assume ($\mathcal{E}$) and (A) hold. Suppose there exists a weak solution $\mu$ of \eqref{SDE}. For each $0 \le s < t \le T$, assume that the function $F_{s,t} : \P(\C) \rightarrow \R$ defined by
\[
F_{s,t}(\nu) := \int_\C\int_s^t\left|\sigma^{-1}b(u,x,\nu) - \sigma^{-1}b(u,x,\mu)\right|^2\,du\,\nu(dx)
\]
is $\tau(\C)$-continuous and $\EE(\P(\C))$-measurable. Assume lastly that there exists $L > 0$ such that
\begin{align}
F_{s,t}(\nu) \le L(t-s)\H(\nu|\mu), \quad \forall s < t, \ \nu \in \P(\C). \label{ass:propagation}
\end{align}
Then the following hold:
\begin{enumerate}
\item For every $\tau(\C)$-open $\EE(\P(\C))$-measurable neighborhood $U$ of $\mu$,
\begin{align}
\limsup_{n\rightarrow\infty}\frac{1}{n}\log \PP(\widehat\mu^n \notin U) \le -e^{-LT}\inf_{\nu \notin U}\H(\nu | \mu). \label{def:LDbound-propagation}
\end{align}
\item For every $\tau(\C)$-open $\EE(\P(\C))$-measurable neighborhood  $U$ of $\mu$, $\lim_{n\rightarrow\infty}\PP(\widehat\mu^n \notin U) = 0$.
\item For each $k \in \N$, we have
\begin{align*}
\lim_{n\rightarrow\infty}\left\|\PP \circ (X^{n,1},\ldots,X^{n,k})^{-1} - \mu^{\otimes k}\right\|_{\mathrm{TV}} = \lim_{n\rightarrow\infty}\H\left( \mu^{\otimes k} \, \big| \, \PP \circ (X^{n,1},\ldots,X^{n,k})^{-1}\right) = 0.
\end{align*}
\end{enumerate}
\end{theorem}

The closest result we know of to Theorem \ref{theorem-propagation} is that of \cite[Theorem 3]{arous1999increasing}, which proves (3) above even when $k$ can grow with $n$, but only when the coefficients (in particular, the interactions) take a very specific form.

The assumption \eqref{ass:propagation} in Theorem \ref{theorem-propagation} is worth commenting on, so we point out two notable sufficient conditions. First, in light of Pinsker's inequality, assumption (B1) is sufficient for \eqref{ass:propagation}.
For a second example, suppose  $\sigma$ is the identity, and the initial law $\lambda_0$ satisfies $\int_{\R^d}\exp(a|x|^2)\lambda_0(dx) < \infty$ for some $a > 0$. Then, by boundedness of $b$ and exponential integrability of Brownian motion, there exists $\tilde{a} > 0$ such that
\[
\int_\C \exp\Big(\tilde a \sup_{t \in [0,T]}|x_t|^2\Big)\,\mu(dx) < \infty.
\]
It follows \cite[Proposition 6.3]{gozlan-leonard} that there exists $C > 0$ such that $\mu$ satisfies the transport inequality
\[
\W_1(\mu,\nu) \le \sqrt{C\H(\nu | \mu)}, \quad \forall \nu \in \P(\C),
\]
where $\W_1$ denotes the $1$-Wasserstein metric on $\P(\C)$.
If we assume $b(t,x,\cdot)$ is Lipschitz with respect to $\W_1$, uniformly in $(t,x)$, then it follows that \eqref{ass:propagation} holds.

\begin{remark}
Conclusion (3) of Theorem \ref{theorem-propagation} implies in particular that
\begin{align}
\lim_{n\rightarrow \infty}\E\left[\prod_{i=1}^k\phi_i(X^{n,i})\right] = \prod_{i=1}^k\int_\C \phi_i \, d\mu, \label{def:classicalpropagation}
\end{align}
for each $k \in \N$ and $\phi_1,\ldots,\phi_k \in B(\C)$. In fact, for fixed $k$, this convergence is uniform over all $\phi_1,\ldots,\phi_k \in B(\C)$ satisfying $|\phi_i| \le 1$.
\end{remark}

\begin{remark} \label{re:convergenceinprob}
Conclusion (2) of Theorem \ref{theorem-propagation} implies that $G(\widehat\mu^n)$ converges to $G(\mu)$ in probability for every $\tau(\C)$-continuous $\EE(\P(\C))$-measurable function $G : \P(\C) \rightarrow \R$. Indeed, for $\epsilon > 0$, $U= \{\nu \in \P(\C) : |G(\nu)-G(\mu)| < \epsilon\}$ is a $\tau(\C)$-open $\EE(\P(\C))$-measurable neighborhood of $\mu$.
\end{remark}

\begin{remark}
The bound (1) is a crude large deviation-type upper bound. The proof employs a change of measure technique reminiscent of the Dawson-G\"artner \cite[Section 5]{dawson-gartner} proof of the large deviation principle for the McKean-Vlasov limit. Following their arguments, one could derive under our same assumptions (even without \eqref{ass:propagation}) the same \emph{local} large deviation bounds as in \cite[Theorem 5.2]{dawson-gartner}, but in the stronger topology $\tau(\C)$. However, to deduce from this a full LDP in the topology $\tau(\C)$ analogous to \cite[Theorem 5.1]{dawson-gartner}, one would need to establish exponential tightness of $\widehat\mu^n$ in the same topology, which does not seem feasible.
\end{remark}

\begin{remark}
It is not true in the setting of Theorem \ref{theorem-propagation} that $\PP(\lim_n \widehat\mu^n = \mu) = 1$, where the limit is taken in $\tau(\C)$. In fact, $\PP(\lim_n \widehat\mu^n = \mu) = 0$, because  for each $\omega \in \Omega$ the countable set $S(\omega) = \{X^{n,i}(\omega) : n \in \N, \ 1 \le i \le n\}$ satisfies both $\widehat\mu^n(\omega)(S(\omega)) = 1$ and $\mu(S(\omega))=0$, as $\mu$ is nonatomic.\footnote{Many thanks to Marcel Nutz for pointing this out.}
In general, a sequence of discrete measures can never $\tau(\C)$-converge to a nonatomic measure, so we cannot hope to improve the convergence in probability stated in Theorem \ref{theorem-propagation}(2). For this reason, we cannot state a version of Theorem \ref{theorem-propagation} in line with more traditional propagation of chaos results (e.g., \cite[Theorem 3.1]{gartner1988mckean}), in which the initial states $X^{n,i}_0$ are taken to be deterministic but with a prescribed limit $\lambda_0 = \lim_n \frac{1}{n}\sum_{k=1}^n\delta_{X^{n,k}_0}$.
\end{remark}

\subsection{A rank-based interaction} \label{se:example}
A notable class of examples related to Burgers' and porous medium type PDEs fits into our framework.
Consider the one-dimensional case $d=1$, with $\sigma \equiv 1$ and
\[
b(t,x,\mu) = g\left(\int_{\C}1_{[0,\infty)}(x_t-y_t)\,\mu(dy)\right),
\]
where $G : [0,1] \rightarrow \R$ is Lipschitz continuous.
The corresponding McKean-Vlasov equation is
\begin{align*}
dX_t = g(\mu_t(-\infty,X_t])dt + dW_t, \quad X_0 \sim \lambda_0, \quad \mu_t = \mathrm{Law}(X_t), \ \forall t \in [0,T].
\end{align*}
Letting $V(t,x) = \mu_t(-\infty,x]$, one expects (cf. \cite{shkolnikov2012large,jourdain1997diffusions,bossy-talay2}) that $V$ is the unique generalized solution of the Burgers-type equation $\partial_t V = \tfrac12\partial_{xx}V - \partial_x(G(V))$, where $G$ is an antiderivative of $g$, and this reduces to Burgers' equation when $g(x)=x$. The corresponding $n$-particle approximation is
\begin{align*}
dX^{n,i}_t = G\left(\frac{1}{n}\sum_{k=1}^n1_{\{X^{n,k}_t \le X^{n,i}_t\}}\right)dt + dW^i_t,
\end{align*}
where $X^{n,i}_0$ are i.i.d.\ with law $\lambda_0$, and $W^i$ are independent Brownian motions. 

All of the assumptions of our Theorems \ref{existence-uniqueness}, \ref{existence}, and \ref{theorem-propagation} hold in this example. Notably,  our Theorem \ref{theorem-propagation}(3) is considerably stronger than \cite[Theorem 3.2]{bossy-talay2} or \cite[Theorem 2.4]{jourdain1997diffusions}, which provide only weak convergence.

\section{Existence and uniqueness proofs} \label{se:proofs-existenceuniqueness}
The proofs of both Theorems \ref{existence-uniqueness} and \ref{existence} rely on the following change of measure argument.  Let $(\Omega,\F,\mathbb{F} = (\F_t)_{0 \le t \le T},P)$ denote a filtered probability space supporting an $\mathbb{F}$-Wiener process $W$ and an $\F_0$-measurable random vector $\xi : \Omega \rightarrow \R^d$ with $P \circ \xi^{-1} = \lambda_0$. Assume $\F_t = \sigma(\xi,W_s : s \le t)$. Let $X$ denote the unique solution of the SDE
\[
dX_t = \sigma(t,X)dW_t, \quad X_0 = \xi.
\]
For each $\mu \in \P(\C)$, define a measure $P_\mu \sim P$ by
\[
\frac{dP_\mu}{dP} := \EE_T\left(\int_0^\cdot \sigma^{-1}b(t,X,\mu) \cdot dW_t\right),
\]
where we define $\EE_t(M) = \exp(M_t - \tfrac12 [M]_t)$ for any continuous martingale $M$.
Girsanov's theorem implies that
\[
W^\mu_t := W_t - \int_0^t\sigma^{-1}b(s,X,\mu)ds
\]
defines a $P_\mu$-Wiener process, and
\[
dX_t = b(t,X,\mu)dt + \sigma(t,X)dW^\mu_t.
\]
Then, a measure $\mu \in \P(\C)$ is a weak solution of \eqref{SDE} if and only if $P_\mu \circ X^{-1} = \mu$.

For $t \in [0,T]$ and $\mu,\nu \in \P(\C)$, abbreviate  $\H_t(\nu | \mu) := \H(\nu^t | \mu^t)$.
Let $\Phi(\mu) := P_\mu \circ X^{-1}$ for $\mu \in \P(\C)$. For any $\mu,\nu \in \P(\C)$, we have
\begin{align*}
\H_t(\Phi(\mu) | \Phi(\nu)) &= -\int_\C \log\frac{d\Phi(\nu)^t}{d\Phi(\mu)^t}d\Phi(\mu)^t = -\E^{P_\mu}\left[ \log\frac{d\Phi(\nu)^t}{d\Phi(\mu)^t}(X_{\cdot \wedge t})\right] \\
	&= -\E^{P_\mu}\left[ \log \E\left[\left.\frac{dP_\nu}{dP_\mu}\right| X_{\cdot \wedge t}\right]\right].
\end{align*}
Assumption (A) and nondegeneracy of $\sigma$ imply that $W$ and $X$ generate the same filtration. Hence, 
\[
\E^{P_\mu}\left[\left.\frac{dP_\nu}{dP_\mu}\right| X_{\cdot \wedge t}\right] = \E^{P_\mu}\left[\left.\frac{dP_\nu}{dP_\mu}\right| \F_t\right],
\]
and so
\begin{align}
\H_t(\Phi(\mu) | \Phi(\nu)) &= -\E^{P_\mu}\left[ \log \E^{P_\mu}\left[\left.\frac{dP_\nu}{dP_\mu}\right| \F_t\right]\right] \nonumber \\
	&= -\E^{P_\mu}\left[ \log \EE_t\left(\int_0^\cdot \left(\sigma^{-1}b(s,X,\nu) - \sigma^{-1}b(s,X,\mu) \right)\cdot dW^\mu_s\right)\right] \nonumber \\
	&= \frac12\E^{P_\mu}\left[\int_0^t\left|\sigma^{-1}b(s,X,\nu) - \sigma^{-1}b(s,X,\mu)\right|^2ds\right]. \label{proof1}
\end{align}

\subsection*{Proof of Theorem \ref{existence-uniqueness}}
We use Banach's fixed point theorem on the complete metric space $(\P(\C),\|\cdot\|_{\mathrm{TV}})$. For any $\mu,\nu \in \P(\C)$, we use \eqref{proof1} along with assumption (B1) to get
\begin{align*}
\H_t(\Phi(\mu) | \Phi(\nu)) &\le \frac12 \kappa^2\int_0^t\|\nu^s - \mu^s\|^2_{\mathrm{TV}}ds.
\end{align*}
By Pinsker's inequality,
\[
\|\Phi(\nu)^t - \Phi(\mu)^t\|^2_{\mathrm{TV}} \le 2\H_t(\Phi(\mu) | \Phi(\nu)) \le \kappa^2\int_0^t\|\nu^s - \mu^s\|^2_{\mathrm{TV}}ds.
\]
Conclude by Picard iteration.\footnote{For the reader worried about measurability of the integrand $s \mapsto \|\nu^s - \mu^s\|^2_{\mathrm{TV}}$, notice that we may write
\[
\|\mu\|_{\mathrm{TV}} = \sup\left\{\int_\C f\,d\mu : f : \C\rightarrow \R \text{ continuous, } |f| \le 1\right\},
\]
from which it is clear that the total variation norm is lower semicontinuous and thus Borel measurable with respect to the topology of weak convergence on $\P(\C)$.}
\hfill\qedsymbol

\subsection*{Proof of Theorem \ref{existence}}
This proof is by Schauder's fixed point theorem, on the topological vector space of bounded signed measures on $\C$ endowed with the weak$^*$ topology induced by $B(\C)$. Note that the induced topology on the subset $\P(\C)$ is exactly $\tau(\C)$.  Proceeding as in \eqref{proof1}, for any $\mu \in \P(\C)$ we have
\begin{align*}
\H(\Phi(\mu) | P \circ X^{-1}) &= \frac12\E^{P_\mu}\left[\int_0^t\left|\sigma^{-1}b(s,X,\mu)\right|^2ds\right] \le \frac12 c^2 T,
\end{align*}
where the constant $c > 0$ comes from assumption (A). Hence,
\[
\Phi(\P(\C)) \subset \left\{\nu \in \P(\C) : \H(\nu | P \circ X^{-1}) \le c^2T/2\right\}.
\]
Sub-level sets of relative entropy are convex, compact, and metrizable in $\tau(\C)$ \cite[Lemma 6.2.12]{dembo-zeitouni}. Hence, to apply Schauder's theorem it remains only to show that $\Phi :\P(\C) \rightarrow \P(\C)$ is sequentially $\tau(\C)$-continuous.
Fix  $\nu,\mu \in \P(\C)$, and use Pinsker's inequality with \eqref{proof1} to get
\begin{align*}
\|\Phi(\nu) - \Phi(\mu)\|^2_{\mathrm{TV}} &\le 2\H(\Phi(\nu) | \Phi(\mu)) = \frac12\int_\C\int_0^T\left|\sigma^{-1}b(s,X,\nu) - \sigma^{-1}b(s,X,\mu)\right|^2\,ds\,\mu(dx).
\end{align*}
As a function of $\nu$, the right-hand side is sequentially $\tau(\C)$-continuous at $\nu=\mu$ by assumption (B2), and this completes the proof.
\hfill\qedsymbol

\section{Proof of Theorem \ref{theorem-propagation}} \label{se:proof-propagation}

We first introduce some notation, used in the proof of both claims (1) and (2).
We transfer the problem set up to a convenient probability space.
Let $(\Omega,\F,P)$ be a probability space supporting an i.i.d. sequence of processes $X^i$ with law $\mu$. For $n \in \N$, let $\mathbb{F}^n = (\F^n_t)_{0 \le t \le T}$ denote the filtration generated by $(X^1,\ldots,X^n)$. There exist i.i.d. Wiener processes $W^1,W^2,\ldots$ such that
\[
dX^i_t = b(t,X^i,\mu)dt + \sigma(t,X^i)dW^i_t,
\]
and such that $W^i$ is adapted to the filtration generated by $X^i$. For $n \in \N$, let
\[
\widehat\mu^n = \frac{1}{n}\sum_{i=1}^n\delta_{X^i}
\]
and
\[
\Delta^i_t := \sigma^{-1}b(t,X^i,\mu^n) - \sigma^{-1}b(t,X^i,\mu).
\]
Define a measure $P^n$ on $(\Omega,\F^n_T)$ by $dP^n/dP = Z^n_T$, where we define the density process
\[
Z^n_t := \EE_t\left(\int_0^\cdot\sum_{i=1}^n\Delta^i_s \cdot dW^i_s \right).
\]
By Girsanov's theorem, $W^{n,i}_\cdot := W^i_\cdot - \int_0^\cdot\sigma^{-1}b(t,X^i,\mu^n)dt$ defines a $P^n$-Wiener process, and
\[
dX^i_t = b(t,X^i,\mu^n)dt + \sigma(t,X^i)dW^{n,i}_t.
\]
Hence $P^n \circ (X^1,\ldots,X^n)^{-1}$ is a weak solution of the $n$-particle system, and in the notation of Section \ref{se:propagation} we have $\PP \circ (X^{n,1},\ldots,X^{n,n})^{-1} = P^n \circ (X^1,\ldots,X^n)^{-1}$.

{\ } \\
\noindent\textbf{Proof of (1).} Fix a $\EE(\P(\C))$-measurable open set $U \subset \P(\C)$ containing $\mu$. The goal is to show that
\begin{align}
\lim_{n\rightarrow\infty}P^n(\mu^n \notin U) = 0. \label{pf:propagation11}
\end{align}
Fix $p,q \in (1,\infty)$, and let $p^*$ and $q^*$ denote the conjugate exponents, $p^*=p/(p-1)$ and $q^*=q/(q-1)$. Assume $p$ and $q$ are such that $M = LTpq/2$ is an integer, for reasons which will be clear later. Define $t_j = jT/M$ for $j=0,\ldots,M$. We will show inductively that, for each $j$,
\begin{align}
\limsup_{n\rightarrow\infty}\frac{1}{n}\log\E^P\left[P^n\left(\mu^n \notin U \, | \, \F_{t_j}\right)\right] \le -(p^*q^*)^{-(M-j)}\inf_{\nu \notin U}\H(\nu | \mu). \label{pf:propagation-induction}
\end{align}
Indeed, once this is established, it is easy to complete the proof of (1) as follows: By taking $j=0$ and noting that $P^n$ and $P$ agree on $\F_{t_0}=\F_0$, it follows from \eqref{pf:propagation-induction} that
\begin{align*}
\limsup_{n\rightarrow\infty}\frac{1}{n}\log P^n(\mu^n \notin U) \le -\left(\frac{pq}{(q-1)(p-1)}\right)^{-LTpq/2}\inf_{\nu \notin U}\H(\nu | \mu).
\end{align*}
Noting that $\lim_{x\rightarrow\infty}(\tfrac{x}{x-1})^x = e$, we may send $p,q \rightarrow \infty$ in the above to get \eqref{def:LDbound-propagation}.

We first check that \eqref{pf:propagation-induction} is valid for $j=M$. By Sanov's theorem \cite[Theorem 6.2.10]{dembo-zeitouni}, $P \circ (\mu^n)^{-1}$ satisfies a large deviation principle (LDP) on $(\P(\C),\tau(\C))$ with good rate function $H(\cdot | \mu)$.
Thus, since $\{\mu^n \notin U\}$ belongs to $\F_{t_M}=\F_T$, we have
\begin{align*}
\limsup_{n\rightarrow\infty}\frac{1}{n}\log\E^P\left[P^n\left(\mu^n \notin U \, | \, \F_{t_M}\right)\right] &= \limsup_{n\rightarrow\infty}\frac{1}{n}\log P(\mu^n \notin U) \\
	&\le -\inf_{\nu \notin U}\H(\nu|\mu).
\end{align*}
We now prove \eqref{pf:propagation-induction} by induction. Suppose \eqref{pf:propagation-induction} holds for some $j=1,\ldots,M$. We then estimate
\begin{align*}
\E^P\left[P^n\left(\mu^n \notin U \, | \, \F_{t_{j-1}}\right)\right] &= \E^P\left[\E^{P^n}\left[\left. P^n\left(\mu^n \notin U \, | \, \F_{t_j}\right) \right| \F_{t_{j-1}}\right]\right] \\
	&= \E^P\left[\E^P\left[\left. \frac{Z^n_{t_j}}{Z^n_{t_{j-1}}} P^n\left(\mu^n \notin U \, | \, \F_{t_j}\right) \right| \F_{t_{j-1}}\right]\right] \\
	&= \E^P\left[\frac{Z^n_{t_j}}{Z^n_{t_{j-1}}} P^n\left(\mu^n \notin U \, | \, \F_{t_j}\right)\right].
\end{align*}
Indeed, the second step follows from Bayes' rule \cite[Lemma 3.5.3]{karatzas-shreve}. Taking note of the identity
\begin{align}
\frac{Z^n_{t_j}}{Z^n_{t_{j-1}}} &= \EE_{t_j}\left(\int_{t_{j-1}}^\cdot\sum_{i=1}^n\Delta^i_s \cdot dW^i_s \right) \nonumber \\
	&= \EE_{t_j}\left(p\int_{t_{j-1}}^\cdot\sum_{i=1}^n\Delta^i_s \cdot dW^i_s \right)^{1/p}\exp\left(\frac{p}{2}\int_{t_{j-1}}^{t_j}\sum_{i=1}^n|\Delta^i_s|^2ds \right)^{1/p^*}, \label{def:RNbound}
\end{align}
we use H\"older's inequality twice to get
\begin{align*}
\E^P&\left[\frac{Z^n_{t_j}}{Z^n_{t_{j-1}}} P^n\left(\mu^n \notin U \, | \, \F_{t_j}\right)\right] \\
	&\le \E^P\left[P^n\left(\mu^n \notin U \, | \, \F_{t_j}\right)^{p^*}\exp\left(\frac{p}{2}\int_{t_{j-1}}^{t_j}\sum_{i=1}^n|\Delta^i_s|^2ds \right)\right]^{\frac{1}{p^*}} \\
	&\le \E^P\left[P^n\left(\mu^n \notin U \, | \, \F_{t_j}\right)^{p^*q^*}\right]^{\frac{1}{p^*q^*}}\E^P\left[\exp\left(\frac{pq}{2}\int_{t_{j-1}}^{t_j}\sum_{i=1}^n|\Delta^i_s|^2ds \right)\right]^{\frac{1}{p^*q}}.
\end{align*}
Hence, 
\begin{align*}
\limsup_{n\rightarrow\infty}\frac{1}{n}\log\E^P\left[P^n\left(\mu^n \notin U \, | \, \F_{t_{j-1}}\right)\right] \le &\frac{1}{p^*q^*}\limsup_{n\rightarrow\infty}\frac{1}{n}\log\E^P\left[P^n\left(\mu^n \notin U \, | \, \F_{t_j}\right)\right]  \\
	&+ \frac{1}{p^*q}\limsup_{n\rightarrow\infty}\frac{1}{n}\log \E^P\left[\exp\left(\frac{pq}{2}\int_{t_{j-1}}^{t_j}\sum_{i=1}^n|\Delta^i_s|^2ds \right)\right].
\end{align*}
In light of the induction hypothesis, the proof will be compete if we show that the last term is not positive.
To do this, we again exploit the fact that $P \circ (\mu^n)^{-1}$ satisfies a LDP on $(\P(\C),\tau(\C))$ with good rate function $H(\cdot \, | \, \mu)$. Using Varadhan's integral lemma
\cite[Theorem 4.3.1]{dembo-zeitouni}  with the $\tau(\C)$-continuous function $F_{t_{j-1},t_j}$,
\begin{align*}
\limsup_{n\rightarrow\infty}\frac{1}{n}\log \E^P\left[\exp\left(\frac{pq}{2}\int_{t_{j-1}}^{t_j}\sum_{i=1}^n|\Delta^i_s|^2ds \right)\right] &= \limsup_{n\rightarrow\infty}\frac{1}{n}\log \E^P\left[\exp\left(\frac{npq}{2}F_{t_{j-1},t_j}(\mu^n) \right)\right] \\
	&\le \sup_{\nu \in \P(\C)}\left(\frac{pq}{2}F_{t_{j-1},t_j}(\nu) - \H(\nu | \mu)\right) \\
	&\le \sup_{\nu \in \P(\C)}\left(\frac{LTpq}{2M} - 1\right)\H(\nu | \mu),
\end{align*}
where we used the assumption \eqref{ass:propagation} and $t_j-t_{j-1}=T/M$.
Recalling that $M = LTpq/2$, the right-hand side equals zero.

{\ } \\
\noindent\textbf{Proof of (2).} Fix a $\EE(\P(\C))$-measurable open set $U \subset \P(\C)$ containing $\mu$. By (1), it suffices to show that $\inf_{\nu \notin U}\H(\nu | \mu) > 0$. But this is a straightforward consequence of the fact that the sub-level set $\{ \nu \in \P(\C) : \H(\nu | \mu) \le a\}$ is $\tau(\C)$-compact for each $a \in \R$ by \cite[Lemma 6.2.16]{dembo-zeitouni}.

{\ } \\
\noindent\textbf{Proof of (3).} 
Define $(\Omega,\F,P)$, $(P^n)$, $(X^i)$, and $(W^{n,i})$ as in the previous step.
Fix $k \in \N$. For $n > k$, define a measure $Q^n$ on $(\Omega,\F^n_T)$ by
\[
\frac{dQ^n}{dP} := \EE_T\left(\int_0^\cdot \sum_{i=k+1}^n\Delta^i_t \cdot dW^i_t\right), \quad \text{ or equivalently }\quad  \frac{dQ^n}{dP^n} = \EE_T\left(-\int_0^\cdot \sum_{i=1}^k\Delta^i_t \cdot dW^{n,i}_t\right).
\]
Note that $Q^n \circ (X^1,\ldots,X^k)^{-1} = P \circ (X^1,\ldots,X^k)^{-1} = \mu^{\otimes k}$. By Pinsker's inequality we have
\begin{align*}
&\left\|P^n \circ (X^1,\ldots,X^k)^{-1} - P \circ (X^1,\ldots,X^k)^{-1}\right\|_{\mathrm{TV}}^2 	 \\
	&\quad\quad\quad\quad  \le 2\H(P \circ (X^1,\ldots,X^k)^{-1} | P^n \circ (X^1,\ldots,X^k)^{-1})  \\
	&\quad\quad\quad\quad  = 2\H(Q^n \circ (X^1,\ldots,X^k)^{-1} | P^n \circ (X^1,\ldots,X^k)^{-1})  \\
	&\quad\quad\quad\quad  = -2\E^{Q^n}\left[\log\E^{Q^n}\left[\left.\frac{dP^n}{dQ^n}\right|X^1,\ldots,X^k\right]\right] \\
	&\quad\quad\quad\quad  = \E^{Q^n}\left[\int_0^T\sum_{i=1}^k\left|\sigma^{-1}b(t,X^i,\mu^n)-\sigma^{-1}b(t,X^i,\mu) \right|^2dt\right] \\
	&\quad\quad\quad\quad  = k\E^{Q^n}\left[\int_0^T\left|\sigma^{-1}b(t,X^1,\mu^n)-\sigma^{-1}b(t,X^1,\mu) \right|^2dt\right].
\end{align*}
Recalling the form of $dQ^n/dP^n$, we use \eqref{def:RNbound} along with Cauchy-Schwarz to bound this by
\begin{align*}
k\E^{P^n}&\left[\exp\left(\int_0^T \sum_{i=1}^k|\Delta^i_t|^2 dt\right)\left(\int_0^T\left|\sigma^{-1}b(t,X^1,\mu^n)-\sigma^{-1}b(t,X^1,\mu) \right|^2dt\right)^2\right]^{1/2} \\
	&\le 4kTc^2e^{4kTc^2}\E^{P^n}\left[\int_0^T\left|\sigma^{-1}b(t,X^1,\mu^n)-\sigma^{-1}b(t,X^1,\mu) \right|^2dt\right]^{1/2},
\end{align*}
where the second step used the bound $|\sigma^{-1}b| \le c$.
Lastly, use symmetry to write
\begin{align*}
\E^{P^n}\left[\int_0^T\left|\sigma^{-1}b(t,X^i,\mu^n)-\sigma^{-1}b(t,X^i,\mu) \right|^2dt\right] &= \E^{P^n}\left[F_{0,T}(\mu^n)\right].
\end{align*}
This converges to zero as $n\rightarrow\infty$ because $F_{0,T}$ is bounded, $\EE(\P(\C))$-measurable, and $\tau(\C)$-continuous; see part (2) and Remark \ref{re:convergenceinprob}.
\hfill\qedsymbol

\bibliographystyle{amsplain}
\bibliography{MKV-bib}

\end{document}